\documentclass{article}

\usepackage{amsfonts}
\usepackage{amsmath}
\mathchardef\mhyphen="2D 
\usepackage[a4paper]{geometry}
\usepackage{microtype}
\usepackage{relsize}

\def\zo/{$0\mkern2mu\mhyphen1$}

\def\nn/{$n \times n$}

\title{Regular Bipartite Lattices with Large Values of $\tht/C_4$}
\date{\today} 
\author{Paul Federbush\\
Department of Mathematics\\
University of Michigan\\
Ann Arbor, MI, 48109-1043}

\usepackage{amsthm} 

\newtheorem*{conj*}{Conjecture} 

\usepackage{tikz}
\newcommand{\tht}{\theta_{2,2,2}}

\begin{document}

\maketitle
\begin{abstract}
The quantities $C_4$ and $\tht$ are as defined by Wanless, $C_4$ just the number of 4-loops of a graph. The construction of this paper provides a counterexample to a conjecture of Butera, Pernici, and the author about the monomer-dimer entropy, $\lambda,$  of a regular bipartite lattice. The lattice we construct is not a lattice graph in its most common definition.
\end{abstract}
We let $G_1 $ and $G_2$ be the graphs illustrated in Figure \ref{fig:graphs}. Then given a graph $G$, we let $C_4= C_4(G)$ be the number of subgraphs of $G$ isomorphic to $G_1,$ and $\tht(G)$ be the number of subgraphs of $G$ isomorphic to $G_2.$ This notation follows Section 2 of \cite{Wan}.

\begin{figure}
\centering
\begin{tikzpicture}
\tikzstyle{v} = [circle, minimum size=2mm,inner sep=0pt,draw]
  \node[v] (n1) at (1.5,2.5) {};
  \node[v] (n2) at (2.5,2.5)  {};
  \node[v] (n3) at (1.5,3.5) {};
  \node[v] (n4) at (2.5,3.5)  {};
  \node (n5) at (2,1.5) {$G_1$};

\foreach \from/\to in {n1/n2,n2/n4,n4/n3,n3/n1}
    \draw (\from) -- (\to);
    
  \node[v] (n6) at (6,2) {};
  \node[v] (n7) at (5,3)  {};
  \node[v] (n8) at (6,3) {};
  \node[v] (n9) at (7,3)  {};
  \node[v] (n10) at (6,4)  {};
  \node (n11) at (6,1) {$G_2$};

\foreach \from/\to in {n6/n7,n6/n8,n6/n9,n7/n10, n8/n10,n9/n10}
    \draw (\from) -- (\to);

\end{tikzpicture}
\caption{}
\label{fig:graphs}
\end{figure}

We propose the following definition for a lattice in the statistical mechanics setting:
A connected infinite graph is a lattice if there is a finitely generated Abelian group of isomorphic maps on $G$ for which the set of equivalence classes of vertices is finite.  Two vertices are in the same equivalence class if there is a translation (a member of the Abelian group) that carries one into the other. The lattice is i-dimensional if exactly i of a set of generators of the group are infinite order.

The lattice ( infinite graph ) we construct in this paper has the following properties:
\vspace{.03in}
\newline1) It is embedded in $R^3$. Edges embed as straight line segments joining vertices.
\newline 2) It is invariant under translations by an integer in the three coordinate directions.
\newline 3) We consider a tessellation of $R^3$ by unit cubes whose vertices have integer coordinates.
Each vertex of the graph is in the interior of one of these cubes. Each cube contains a finite number of vertices.
\newline 4) Each edge either lies in the interior of some cube, or connects two vertices in nearest neighbor cubes.
\vspace{.1in}
\newline This lattice does not have the property of arising from a tessellation, and thus is not a lattice graph
as most commonly understood, though I believe it fits the idea of being a lattice for many researchers.

For a lattice the quantities $C_4$ and $\tht$ may both be infinite, but the ratio $\tht/C_4$ makes sense. Taking a sequence of finite graphs that approach the lattice in the sense of Benjamini-Schramm \cite{CAH}
the ratios $\tht/C_4$ converge to a limit we take to be the value for the lattice.

Given any number $\kappa$ we find a three dimensional regular bipartite lattice for which \begin{align} \tht/C_4 > \kappa \\ C_6 = 0.\end{align} Here $C_6$ is the number of 6-loops as $C_4$ is the number of 4-loops.

In \cite{FF} Friedland and I introduced an expansion for the monomer-dimer entropy of a regular lattice, \begin{align}
\lambda(p)=\lambda_B(p) + \sum_k d_k p^k\end{align} where $\lambda_B$ is the monomer-dimer entropy for the Bethe lattice. In \cite{BFP} Butera, Pernici and I conjectured that for a regular bipartite lattice all the $d_k\ge 0$. Actually a slightly weaker
conjecture was stated in \cite{BFP}, see \cite{Csi} Section 3.
In \cite{Per} Pernici proved in the bipartite case $d_2,d_3,d_4,d_5$ were all nonnegative using a formalism of Wanless, \cite{Wan}. For $d_6$ he finds the expression \begin{align}
d_6 = \frac{5 \bar{C_4}}{d^6} + \frac{\bar{C_6}}{2d^6} - 2\frac{\bar\theta_{2,2,2}}{d^6}
\end{align}
where we have d-regularity and the bars indicate an average over the number of vertices, again as may be taken as a Benjamini-Schramm limit. We find a three dimensional regular bipartite lattice for each $\kappa$ such that $$\tht/C_4 > \kappa$$ and also $$C_6=0.$$ Thus there are bipartite regular lattices for which $d_6 < 0.$ Pernici erred in stating otherwise.

However the conjecture that all the $d_i \ge 0$ should not be lightly dismissed. For hyper-cubic lattices of all dimensions it is
shown in \cite{BFP} that the
first $20$ such coefficients satisfy this condition! It seems quite possible there is a class of 'nice' lattices that satisfy the
conjecture, possibly vertex-transitive lattices.

We start the construction of our lattice with a `root unit graph', a finite graph. To this we apply a number of 2-lifts to end with a `full unit graph'. The full unit graphs are attached to yield the three dimensional lattice. We use 2-lifts similarly to as in Section 4 of \cite{Cs2}, derived from \cite{8}. The construction requires that $d \ge 5$.

\vspace{.2 in}

\noindent\underline{The Root Unit Graph}

\vspace{.1 in}

This graph has black vertices: \begin{align*} c_i, &\quad i=1,\ldots,d\end{align*} and white vertices: \begin{align*} t,b,lx,ly,lz,
rx,ry,rz, & \\ f_i,&\quad i=1,\ldots,d-5 \end{align*}
The edges are: \begin{align*}
(t,c_i)& \quad i=1,\ldots, d \\
(b,c_i)&  \quad i=1,\ldots, d \\
(c_i,f_j)& \quad i=1,\ldots, d\quad j=1,\ldots,d-5 \\
&(lx,c_1), (ly,c_1), (lz,c_1) \\
(rx,c_i)& \quad i=2\ldots, d \\
(ry,c_i)& \quad i=2,\ldots, d \\
(rz,c_i)& \quad i=2,\ldots, d \\
\end{align*}

It has a subgraph, the `root central subgraph'. The root central subgraph has black vertices: \begin{align*}
 c_i,& \quad i=1,\ldots,d 
\end{align*} and white vertices: $$t,b$$ Its edges are:
\begin{align*}
(t,c_i)&\quad i=1,\ldots, d \\
(b,c_i)&\quad  i=1,\ldots, d \\
\end{align*}
The quantities $\tht$ and $C_4$ for the root central subgraph are easily computed to be 
\begin{align}
C_4 &= \frac{d(d-1)}{2} \\
\tht &= \frac{d(d-1)(d-2)}{6}
\end{align}
Exploitation of the root central subgraph is the key idea in this paper.

\vspace{.2 in}
\noindent\underline{Thank Heaven for 2-lifts}
\vspace{.1 in}

The full unit graph is constructed from the root unit graph through a sequence of $s$  2-lifts. Each edge in the full unit graph projects in a natural sense onto an edge of the root unit graph. $2^s$ edges of the full unit graph project onto each edge of the root unit graph (likewise $2^s$ vertices project onto each vertex). The edges that project onto edges of the root central subgraph have all been chosen parallel, none crossed (this is also in a natural sense). Thus there are $2^s$ disjoint subgraphs of the full unit graph isomorphic to the root central subgraph. After $s$ 2-lifts any vertex may be assigned a binary integer of length $s$, that we call its 'level'. If a zero appears in the
$i^{th}$ spot then at the $i^{th}$ 2-lift the lower vertex was chosen, if a one appears the upper vertex chosen. ( One arrives at the vertex
by following a sequence of $s$ splittings. )

The sequence of 2-lifts are chosen so that for the full unit graph there are no 6-loops, and 4-loops only in the $2^s$ subgraphs isomorphic to the root central subgraph.

\vspace{.2 in}
\noindent\underline{Connecting the Graph}
\vspace{.1 in}

We associate a full unit graph to each of the aforementioned unit cubes that tessellate $R^3$. Corresponding to each nearest
neighbor pair of cubes we create a connection. Let cubes $a$ and $b$ be nearest neighbors, and suppose cube $b$ is
obtained from cube $a$ by the translation $x \rightarrow x+1$. Then we connect the full unit graph associated to $a$ to the
full unit graph associated to $b$  by identifying each of the $2^s$ vertices that project onto vertex $rx$ of the  graph $a$ with one of the $2^s$ vertices that project onto $lx$ of the graph $b$ in a level preserving way.  Nearest neighbor cubes lined
up in the other two directions are treated analagously. The lattice or graph thus constructed will have $C_6=0$ and  $\frac{\bar \theta _{2,2,2}}{\bar C_4} = \frac{g-2}{3}.$

\vspace{.2 in}
\noindent\underline{Locating the Graph}
\vspace{.1 in}

It remains to associate to each vertex of the lattice a point in $R^3$. To do this it is clearly sufficient to locate each vertex of
a single full unit graph included in the construction. We let $G$ be the full unit graph associated to the unit cube
$[0,1] \times [0,1] \times [0,1]$, and $g$ its root unit graph. We divide the vertices of $g$ into five disjoint sets
$S_1$, $S_2$, $S_3$, $S_4$ and $S_5$with $S_1 = \left \{ rx \right \}$,  $S_2 = \left \{ ry \right \}$ , $S_3 = \left \{ rz \right \}$ and  $S_4 = \left \{ lx,ly,lz \right \}$. We let $\bar {S_i}$
be the vertices in $G$ that project to vertices in $S_i$. A 'try' is an assignment of the vertices in $\bar{S_5}$ into a set of disjoint
points in $(1/3,2/3) \times (1/3,2/3) \times (1/3,2/3)$, an assignment of $\bar{S_1}$ into a set of disjoint points in
  $(2/3,1) \times (1/3,2/3) \times (1/3,2/3)$, an assignment of $\bar{S_2}$ into a set of disjoint points in
  $(1/3,2/3) \times (2/3,1) \times (1/3,2/3)$, and an assignment of $\bar{S_3}$ into a set of disjoint points in
  $(1/3,2/3) \times (1/3,2/3) \times (2/3,1)$. Note that the assignments for $\bar{S_1}$, $\bar{S_2}$, and $\bar{S_3}$
  determine an assignment for $\bar{S_4}$ by the way connections were determined. For example, if a vertex in $\bar{S_1}$
  is assigned a point $\bold w$ in $R^3$, then the vertex of the same level that projects onto $lx$ is assigned the
  point $\bold w - (1,0,0)$. A try is a "good try' if the edges of our full unit graph map into line segments that may intersect
  only in a vertex. But almost all tries are good tries, we pick a good try to fix the assignment. The resulting 
  embedding of the lattice satisfies all our conditions.

\textbf{Acknowledgement}  We would like to thank Sasha Barvinok for a lecture on 2-lifts and John Stembridge for putting up 
with my continual questions.

\end{document}